\documentclass[a4paper,oneside]{amsart}
\usepackage{graphicx}
\usepackage{amssymb}
\usepackage[T1]{fontenc}
\usepackage[latin1]{inputenc}

\makeatletter
 \theoremstyle{plain}
\newtheorem{thm}{Theorem}[section]
  \theoremstyle{plain}
  \newtheorem{prop}[thm]{Proposition}
  \theoremstyle{plain}
  \newtheorem{lem}[thm]{Lemma}
 \theoremstyle{definition}
  \newtheorem{example}[thm]{Example}

\usepackage{a4wide}

\makeatother
\begin{document}
\newcommand{\cat}[1]{\mathsf{#1}}

\newcommand{\map}{\rightarrow}

\newcommand{\Hom}{\mathrm{Hom}}

\newcommand{\HOM}{\mathbf{Hom}}

\newcommand{\rodd}[1]{\mathbb{R}^{0|#1}}

\title[$L_{\infty}$-algebras as 1-jets of simplicial manifolds]{$L_{\infty}$-algebras as 1-jets of simplicial manifolds\\
 (and a bit beyond)}

\author{Pavol \v Severa}

\curraddr{Dept.~of Theoretical Physics, Comenius University, Mlynská dolina
F2, 84248 Bratislava, Slovakia}

\maketitle
\tableofcontents{}

\section{Introduction}

The Lie algebra of a Lie group can be seen, in certain sense, as its
first jet. As a generalization, we shall define first and also higher
jets of arbitrary simplicial manifolds. The first jet will be, under
certain representability condition, an $L_{\infty}$-algebra. The
condition is satisfied, in particular, for Kan simplicial manifolds.
While higher-order jets don't give anything new in the case of a Lie
group, for general simplicial manifolds they do.

Any simplicial manifold (contravariant functor from the category of
finite ordered sets to the category of manifolds) gives rise to a
contravariant functor \[
\mbox{category of surjective submersions}\map\mbox{category of sets}.\]
The idea is to restrict such a functor to submersions with fibres
isomorphic to $\rodd{1}$. As we shall see, under a representability
condition the restricted functor is equivalent to a differential graded
manifold, i.e.~(more-or-less) to an $L_{\infty}$-algebra. We shall
also consider $\rodd{n}$ in place of $\rodd{1}$ (yielding the {}``higher
jets'' mentioned above). The whole {}``theory'' is just an elementary
and simple abstract nonsense - restricting a functor to a full subcategory
and then inducing it back.

The ideas of this paper were already informally presented in \cite{ks}.
The basic idea, that differential forms on $M$ are functions on the
space of maps from $\rodd{1}$ to $M$, and that the structure of
a complex on $\Omega(M)$ is equivalent to the action of $\mathit{Diff}(\rodd{1})$
on this map space, is taken from \cite{kon}.

\section{Introductory example - Lie groups and Lie algebras\label{sec:Introductory-example}}

Let $G$ be a Lie group and $M\map N$ a surjective submersion. Recall
that a $G$-descent data (i.e.~a descent of the trivial principal
$G$-bundle over $M$ to a principal $G$-bundle over $N$) is a map\[
g:M\times_{N}M\map G\]
such that $g(x,x)=e$ and $g(x,y)g(y,z)=g(x,z)$ for any $(x,y,z)\in M\times_{N}M\times_{N}M$.

Let \[
\Omega(M\map N)\]
 denote fibrewise differential forms on $M$. A $G$-descent data
gives us a flat fibrewise $\mathfrak{g}$-connection (where $\mathfrak{g}$
is the Lie algebra of $G$), i.e.~an $\alpha\in\Omega^{1}(M\map N)\otimes\mathfrak{g}$
satisfying the Maurer-Cartan equation\[
d\alpha+[\alpha,\alpha]/2=0.\]
Vice versa, if $G$ and all the fibres are 1-connected, a fibrewise
flat connection gives us descent data.

Let us now describe how the method used in this paper works in this
case; details are provided in the subsequent sections. Let $F(M\map N)$
denote the set of all $G$-descent data on $M\map N$. $F$ is a contravariant
functor from the category of surjective submersions to the category
of sets. Let us restrict $F$ to the submersions of the form $\rodd{1}\times N\map N$
(where the arrow is the canonical projection to the second factor)
and denote it $\mathrm{res}_{1}F$. This restricted functor is equivalent
to the Lie algebra $\mathfrak{g}$: 

\begin{enumerate}
\item as a functor in $N$, $\mathrm{res}_{1}F$ is representable, namely
by the supermanifold $\Pi\mathfrak{g}$ (recall that by definition,
$C^{\infty}(\Pi\mathfrak{g})=\bigwedge\mathfrak{g}^{*}$)
\item the functoriality of $\mathrm{res}_{1}F$ with respect to all commutative
squares\[
\begin{array}{ccc}
\rodd{1}\times N_{1} & \map & \rodd{1}\times N_{2}\\
\downarrow &  & \downarrow\\
N_{1} & \map & N_{2}\end{array}\]
is equivalent, as we shall see, to the structure of a non-negatively
graded differential manifold on $\Pi\mathfrak{g}$, i.e.~to the Lie
algebra structure on $\mathfrak{g}$.
\end{enumerate}
If we induce (where by induction we mean the right adjoint to the
restriction) $\mathrm{res}_{1}F$ back to a functor defined on all
surjective submersions, it will send a submersion to the set of all
flat fibrewise connections; the map\[
G\mbox{-descent data}\mapsto\mbox{flat fibrewise }\mathfrak{g}\mbox{-connection},\]
together with its universal property, then comes just from the fact
that restriction and induction are mutually adjoint.

Generally, if $F$ is any contravariant functor from surjective submersions
to sets, we can apply the same procedure. If $\mathrm{res}_{1}F$
is representable as a functor in $N$ by a supermanifold $X$ then
$X$ becomes a differential non-negatively graded manifold. If $X$
is in fact positively graded, it can be equivalently described as
a non-positively graded finite-dimensional $L_{\infty}$-algebra.
In this case, any element of $F(M\map N)$ would give us a {}``flat
fibrewise $X$-connection'', i.e.~an $\alpha\in\left(\Omega(M\map N)\otimes X\right)^{1}$
such that\[
d\alpha+[\alpha,\alpha]/2!+[\alpha,\alpha,\alpha]/3!+\dots=0.\]
For general $X$ ({}``$L_{\infty}$-algebroid'') we would get a
DGA morphism \[
C^{\infty}(X)\map\Omega(M\map N)\]
 (this is of course equivalent to the above {}``flatness'' condition
in the case of positively graded $X$) for every element of $F(M\map N)$.

An example of such a functor is given by any simplicial manifold $K_{\bullet}$:
For any surjective submersion $M\map N$ let $(M\map N)_{\bullet}$
denote the nerve of the groupoid $M\times_{N}M\rightrightarrows M$.
Let us  consider simplicial morphisms\[
(M\map N)_{\bullet}\map K_{\bullet};\]
 if $K_{\bullet}$ is the nerve of a Lie group $G$, these are just
$G$-descent data on $M\map N$. Let $F(M\map N)$ be the set of these
morphisms. If $K_{\bullet}$ satisfies Kan condition then $\mathrm{res}_{1}F$
is represenentable as a functor in $N$ and we get a differential
graded manifold (a non-positively graded $L_{\infty}$-algebra if
$K_{0}=\{ pt\}$).

\section{Notation}

If $M\map N$ is a surjective submersion, we denote\[
\Omega(M\map N)\]
the space of fibrewise differential forms on $M$, and we shall use
similar notation for other fibrewise objects:\[
\mathcal{Z}^{k}(M\map N)\]
will be the space of closed fibrewise $k$-forms,\[
T(M\map N)\]
will be the fibrewise tangent bundle (a subbundle of $TM$) etc.

\section{Reminder on presheaves (generalized objects)}

The basic reference for this section is \cite{sga4}, exposé I. Let
$\mathsf{C}$ be a category. A \emph{presheaf} on $\cat{C}$ is a
functor $\cat{C}^{o}\map\cat{Set}$; the category of presheaves%
\footnote{we make the usual hyper-correct assumption of working in some universe
of sets to avoid set-of-all-sets-like problems%
} on $\cat{C}$ (with natural transformations as morphisms) is denoted
$\hat{\cat{C}}$. Presheaves can be reasonably viewed as generalized
objects of $\cat{C}$, with $F(X)$ ($F\in\hat{\cat{C}}$, $X\in\cat{C}$)
interpreted as the set of morphisms $X\map F$. Namely, any object
$Y\in\cat{C}$ gives us a presheaf $Y\in\hat{\cat{C}}$ via $Y(X)=\Hom_{\cat{C}}(X,Y)$.
For any $X\in\cat{C}$ and $F\in\hat{\cat{C}}$ then $\Hom_{\hat{\cat{C}}}(X,F)\cong F(X)$.
This way $\cat{C}$ is identified with a full subcategory of $\hat{\cat{C}}$
(which is the excuse for denoting $X\in\cat{C}$ and the corresponding
presheaf $X\in\hat{\cat{C}}$ by the same letter). 

A presheaf isomorphic to some $X\in\cat{C}$ is said to be \emph{representable}.
For example, if $U$, $V$ are objects of $\cat{C}$, the presheaf
$U\times V$ is defined as $U\times V(X)=U(X)\times V(X)$; if it
is representable, the corresponding object of $\cat{C}$ (defined
up to a unique isomorphism) is called the cartesian product of $U$
and $V$, and denoted (somewhat abusively) $U\times V$ as well. Similarly,
the presheaf $\HOM(U,V)$ is defined by $\HOM(U,V)(X)=\Hom(U\times X,V)$.
If it is representable, the corresponding object of $\cat{C}$ is
called the internal Hom from $U$ to $V$, and is still denoted $\HOM(U,V)$.

\section{Differential forms as functions on map spaces\label{sec:forms-worms}}

In this section we work in the category $\cat{SM}$ of all smooth
supermanifolds. We shall first notice that $\HOM(\rodd{1},X)$ is
representable for any $X\in\cat{SM}$, and in fact is isomorphic to
something well known:

\begin{prop}
\label{pro:basic}For any supermanifold $X$, $\HOM(\rodd{1},X)$
is naturally isomorphic with $\Pi TX$, and thus (using induction)
$\HOM(\rodd{k},X)$ is naturally isomorphic with $(\Pi T)^{k}X$.
\end{prop}
Here $\Pi TX$ is the odd tangent bundle of $X$. Recall that differential
forms on $X$ are functions on $\Pi TX$.%
\footnote{more precisely, functions that are polynomial on the fibres of $\Pi TX$;
general functions on $\Pi TX$ are called pseudodifferential forms
on $X$. If $X$ is a manifold then every pseudodifferential form
is actually a differential form.%
} The proof of the proposition is straightforward: let $\theta$ be
the coordinate on $\rodd{1}$ and $x^{i}$ be local coordinates on
$X$, so that $x^{i}$, $dx^{i}$ are local coordinates on $\Pi TX$.
If $\rodd{1}\map X$ is a map parametrized by some $Y$ (that is,
a map $\rodd{1}\times Y\map X$), we expand it to Taylor series in
$\theta$, \[
x^{i}(\theta)=x^{i}(0)+\xi^{i}\theta;\]
 identifying $x^{i}(0)$ with $x^{i}$ and $\xi^{i}$ with $dx^{i}$
we get the (local) isomorphism between $\HOM(\rodd{1},X)$ and $\Pi TX$.
This identification is independent of the choice of coordinates $x^{i}$,
since using other coordinates $\tilde{x}^{i}$ we have \[
\tilde{x}(\theta)=\tilde{x}(x(\theta))=\tilde{x}(x(0)+\xi\theta)=\tilde{x}(x(0))+\frac{\partial\tilde{x}}{\partial x}\xi\theta,\]
i.e.~$\tilde{\xi}=\frac{\partial\tilde{x}}{\partial x}\xi$, just
as $d\tilde{x}=\frac{\partial\tilde{x}}{\partial x}dx$.

The isomorphism of $\HOM(\rodd{1},X)$ with $\Pi TX$ is certainly
not surprising: the relation $\theta^{2}=0$ basically says that maps
$\rodd{1}\map X$ are 1-jets of curves in $X$, i.e.~tangent vectors
in $X$. It gives us, however, an interesting explanation of the de
Rham differential (which, being a derivation on $\Omega(X)$, is a
vector field on $\Pi TX$). Namely, on $\HOM(\rodd{1},X)$ we have
a right action of the supersemigroup $\HOM(\rodd{1},\rodd{1})$, and
this action gives us the structure of a complex on $\Omega(X)$ (a
left action of $\HOM(\rodd{1},\rodd{1})$ on a vector space $V$ is
the same as a structure of a non-negatively graded complex on $V$).

Let us compute the action of $\HOM(\rodd{1},\rodd{1})$ explicitly.
Given a transformation $\theta\mapsto\theta'=a\theta+\beta$ of $\rodd{1}$,
we have\[
x^{i}(\theta')=x^{i}(a\theta+\beta)=\left(x^{i}(0)+\xi^{i}\beta\right)+a\xi^{i}\theta,\]
i.e.\[
x^{i}\mapsto x^{i}+dx^{i}\,\beta,\qquad dx^{i}\mapsto a\, dx^{i}.\]
If we take infinitesimal generators of $\HOM(\rodd{1},\rodd{1})$,
$\partial/\partial\theta$ and $\theta\partial/\partial\theta$, we
get the following:

\begin{prop}
The vector fields $\partial/\partial\theta$ and $\theta\,\partial/\partial\theta$
on $\rodd{1}$ act on $\Pi TM$ as\[
dx^{i}\frac{\partial}{\partial x^{i}}\quad\textrm{and}\quad dx^{i}\frac{\partial}{\partial(dx^{i})},\]
i.e.~as the de Rham differential and the degree.
\end{prop}
In general, a supermanifold $X$ with a right action of $\HOM(\rodd{1},\rodd{1})$
is a \emph{differential non-negatively graded supermanifold}, i.e.
the algebra of functions on $X$ is a differential non-negatively
graded superalgebra.

In Proposition \ref{pro:basic} we noticed that for any $k\in\mathbb{N}$
the iterated odd tangent bundle $(\Pi T)^{k}X$ is isomorphic to $\HOM(\rodd{k},X);$
as a result, we have a right action of the supersemigroup\[
\HOM(\rodd{k},\rodd{k})\]
on $(\Pi T)^{k}X$. Functions on $(\Pi T)^{k}X$ were called \emph{differential
worms} in \cite{ks} and denoted $\Omega_{[k]}(X)$; they were later
introduced under a saner name \emph{iterated differential forms} in
\cite{vv}, from a somewhat different perspective.

\section{Degrees and parity}

This section is just a technical remark on the relation between degrees
and parity.

Let $V$ be a $\mathbb{Z}$-graded vector space. One usually endows
$V$ with a $\mathbb{Z}/2\mathbb{Z}$-grading (parity) using the parity
of the degrees. If we instead suppose that each $V_{k}$ is $\mathbb{Z}/2\mathbb{Z}$-graded,
we shall call $V$ a \emph{$\mathbb{Z}$-graded vector superspace}.

In a graded commutative algebra $A$ we have $ab=(-1)^{\deg a\deg b}ba$.
In a graded commutative \emph{superalgebra} the sign is given by the
parities, rather than by the degrees. In a differential graded commutative
superalgebra we need to say that the differential is both odd and
degree-1.

Let us now look at differential non-negatively graded supermanifolds,
i.e.~at supermanifolds with a right action of $\HOM(\rodd{1},\rodd{1})$.
The algebra of functions on such a supermanifold $X$ is a differential
graded commutative superalgebra. If it is a graded commutative algebra,
i.e.~if the parity is given by the degrees, we shall call $X$ (as
is customary) a \emph{differential non-negatively graded manifold}.
It happens iff the parity involution of $\rodd{1}$ (a map $\rodd{1}\map\rodd{1}$)
acts on $X$ as the parity involution of $X$.

\section{$L_{\infty}$-algebras and differential graded manifolds}

Recall that an $L_{\infty}$-structure on a graded vector space $V$
is a derivation $Q$ of the coalgebra $S(V[1])$ of degree $-1$,
such that $Q^{2}=0$ and $Q1=0$. If we compose $Q$ with the projection
$S(V[1])\map V[1]$, and take the homogeneous components of the resulting
map, we get linear maps\begin{eqnarray*}
Q_{1} & : & V\map V[1]\\
Q_{2} & : & V\wedge V\map V\\
Q_{3} & : & V\wedge V\wedge V\map V[-1]\\
 & \cdots\end{eqnarray*}
The derivation $Q$ is uniquely determined by these components (they
need to satisfy some quadratic equations so that $Q^{2}=0$). $Q_{1}$
makes $V$ into a complex, and is usually denoted $d$; $Q_{2}(u,v)$
is denoted $[u,v]$, $Q_{3}(u,v,w)$ is denoted $[u,v,w]$ etc. If
all $Q_{k}$'s vanish for $k\geq3$ then $d$ and $[,]$ make $V$
into a differential graded Lie algebra (DGLA).

This definition can be rephrased in more geometrical terms (see e.g.~\cite{kon}):
$Q$ is a vector field (of degree 1, vanishing at the origin and such
that $Q^{2}=0$) on the formal graded manifold $V[1]$. {}``Formal''
means that $Q$ is given by a formal power series; $Q_{k}$'s are
its homogeneous components.

In our case $V$ will always be non-positively graded and finite-dimensional.
For dimensional reasons only finitely many $Q_{k}$'s can be non-zero,
so that $Q$ is a vector field on the true (non-formal) graded manifold
$V[1]$. Moreover, any positively graded manifold is of the form $V[1]$
for some non-positively graded finite-dimensional vector space $V$;
$L_{\infty}$-algebras of our type are thus equivalent to differential
positively graded manifolds.

The differential graded manifolds that we shall construct will be,
in general, just non-negatively graded (i.e.~not necessarily positively
graded). Since they are true (non-formal) supermanifolds, they are
not equivalent to $L_{\infty}$-algebras (though $L_{\infty}$-algebras
appear as formal neighbourhoods of invariant points). They might be
called $L_{\infty}$-algebroids.

This section remains valid if we substitute $L_{\infty}$-algebras
with $L_{\infty}$-superalgebras and graded manifolds with graded
supermanifolds.

\section{Simplicial manifolds and presheaves on surjective submersions}

Let $\cat{SSM}$ denote the category of surjective submersions, i.e.~its
objects are surjective submersions between \emph{supermanifolds} and
morphisms are commutative squares.

Any surjective submersion $M\map N$ gives us a simplicial (super)manifold
$(M\map N)_{\bullet}$, namely the nerve of the groupoid $M\times_{N}M\rightrightarrows M$.
Any simplicial (super)manifold $K_{\bullet}$ therefore gives us a
presheaf on $\cat{SSM}$, defined by\[
F_{K_{\bullet}}(M\map N)=\Hom((M\map N)_{\bullet},K_{\bullet}).\]

As a technical nonsense, we shall call a presheaf on $\cat{SSM}$
\emph{even} if it maps the parity involution to the identity. A representable
presheaf is even iff in the coresponding object $M\map N$ both $M$
and $N$ are manifolds. If $K_{\bullet}$ is a simplicial manifold
then $F_{K_{\bullet}}$ is even. The reason for this definition is
that even presheaves will lead to differential graded manifolds ($L_{\infty}$-algebras)
rather than to more general differential graded supermanifolds ($L_{\infty}$-superalgebras).

\section{Jets of presheaves\label{sec:Jets}}

Let $\cat{SSM}_{n}$ denote the full subcategory of $\cat{SSM}$ with
objects $\rodd{n}\times N\map N$ ($N$ runs through all supermanifolds,
the arrow is the canonical projection).

We have \[
\Hom\left(\rodd{n}\times N_{1}\map N_{1},\rodd{n}\times N_{2}\map N_{2}\right)\simeq\Hom(N_{1},N_{2})×\HOM(\rodd{n},\rodd{n})(N_{1}),\]
 so we get the following lemma:

\begin{lem}
\label{lem:PSMn}An object of $\widehat{\cat{SSM}_{n}}$ can be equivalently
described as an object of $\widehat{\cat{SM}}$ with a right action
of $\HOM(\rodd{n},\rodd{n})$. The two categories (\/$\widehat{\cat{SSM}_{n}}$
and the right-$\HOM(\rodd{n},\rodd{n})$-equivariant version of $\widehat{\cat{SM}}$)
are equivalent.
\end{lem}
Let us denote the right-$\HOM(\rodd{n},\rodd{n})$-equivariant version
of $\widehat{\cat{SM}}$ by $\widehat{\cat{SM}}_{[n]}$ and the right-$\HOM(\rodd{n},\rodd{n})$-equivariant
version of $\cat{SM}$ by $\cat{SM}_{[n]}$ ($\cat{SM}_{[n]}$ is
thus a full subcategory of $\widehat{\cat{SM}}_{[n]}$ and $\widehat{\cat{SM}}_{[n]}$
is equivalent to $\widehat{\cat{SSM}_{n}}$). In particular, $\cat{SM}_{[1]}$
is the category of differential graded supermanifolds.

If $F$ is any presheaf on $\cat{SSM}$, we can restrict it to $\cat{SSM}_{n}$,
and we shall call the corresponding object of $\widehat{\cat{SM}}_{[n]}$
the \emph{$n$-jet of $F$.} 

Especially interesting is the case when the $n$-jet of $F$ is representable
(as a presheaf on $\cat{SM}$), when it gives us a supermanifold with
a right $\HOM(\rodd{n},\rodd{n})$-action (one can easily see that
if the $n$-jet of $F$ is representable, so are all its $m$-jets
for $m\leq n$). In the case of $n=1$ we get a differential graded
supermanifold; if $F$ is even, we get a differential graded manifold.

\begin{prop}
\label{pro:Kan}Suppose $K_{\bullet}$ is a simplicial manifold satisfying
Kan condition, which is moreover $m$-truncated for some $m\in\mathbb{N}$,
and denote $F_{K_{\bullet}}$ the corresponding presheaf on $\cat{SSM}$.
Then for every $n\in\mathbb{N}$ the $n$-jet of $F_{K_{\bullet}}$
is representable.
\end{prop}
The proof is a bit technical and can be found (together with definitions)
in Appendix \ref{sec:app}. Let us here just remark that the supermanifold
representing the $n$-jet of $F_{K_{\bullet}}$ is constructed as
a subsupermanifold of $\HOM((\rodd{n})^{m},K_{m})$. In particular
in the case of $n=1$ it is a graded submanifold of $T[1]^{m}K_{m}$
(the differential, however, is not inherited from $T[1]^{m}K_{m}$);
for example, if $K_{\bullet}$ is the nerve of a Lie group $G$, we
get $\mathfrak{g}[1]=T_{e}[1]G\subset T[1]G$ and in this case $K_{1}=G$.
The highest degree of a coordinate is thus at most $m$; the corresponding
$L_{\infty}$-algebra has all degrees higher that $-m$.

If the $n$-jet of an $F\in\widehat{\cat{SSM}}$ is representable
by some $X\in\cat{SM}_{[n]}$, we shall define $\Omega_{[n]}(F)$
to be\[
\Omega_{[n]}(F):=C^{\infty}(X);\]
the reason is that if $F$ itself is representable by $M\map N$ then
this definition gives $\Omega_{[n]}(F)=\Omega_{[n]}(M\map N)$. This
follows from Prop.~\ref{pro:basic} if $N$ is a point, or from its
obvious fibrewise version in general.

\section{Reminder on presheaves (continued)}

Let $u:\cat{C}\map\cat{D}$ be any functor. It induces a functor $u^{*}:\hat{\cat{D}}\map\hat{\cat{C}}$
via $u^{*}(F)=F\circ u$. The functor $u^{*}$ admits a right adjoint
$u_{*}:\hat{\cat{C}}\map\hat{\cat{D}}$, which can be defined via
$u_{*}(F)(Z)=\Hom_{\hat{\cat{C}}}(u^{*}(Z),F)$ where $F\in\hat{\cat{C}}$
and $Z\in\cat{D}$. If now $\cat{C}$ is a full subcategory%
\footnote{There is, of course, no real difference between a fully faithful functor
and a full subcategory. Sometimes we shall even commit the crime of
calling a category $\cat{C}$ to be a full subcategory of $\cat{D}$
when all we have is a fully faithful functor $\cat{C}\map\cat{D}$,
provided the functor is clear from the context ($\cat{C}$ is then
just equivalent to a full subcategory of $\cat{D}$).%
} of $\cat{D}$ and $u$ is the inclusion (the only situation we shall
meet; in that case $u^{*}$ is simply the restriction) then $u^{*}\circ u_{*}$
is (isomorphic to) $\mathrm{id}_{\hat{\cat{C}}}$ and consequently
$u_{*}$ is fully faithful. We can thus use $u_{*}$ to identify $\hat{\cat{C}}$
with a full subcategory of $\hat{\cat{D}}$, and $u^{*}$ gives us
a projection $\hat{\cat{D}}\map\hat{\cat{C}}$.

We finish with a simple condition that forces $u^{*}$ (and thus $u_{*}$)
to be an equivalence of categories.

\begin{lem}
\label{lem:trivequiv}Suppose that $u:\cat{C}\map\cat{D}$ is a fully
faithful functor, and that any object $Y\in\cat{D}$ is a retract
of some object $X\in\cat{C}$, i.e.~that there are morphisms $Y\map u(X)\map Y$
that compose to $\mathrm{id}_{Y}$. Then $u^{*}$, and consequently
$u_{*}$, is an equivalence of categories.
\end{lem}

\section{Approximations and representability of presheaves}

Let $u_{n}:\cat{SSM}_{n}\map\cat{SSM}$ be the inclusion (see Section
\ref{sec:Jets}). Recall that for any presheaf $F$ on $\cat{SSM}$
its restriction $u_{n}^{*}F$ (a presheaf on $\cat{SSM}_{n}$) is
equivalent to a right-$\HOM(\rodd{n},\rodd{n})$-equivariant presheaf
on $\cat{SM}$, which we called the $n$-jet of $F$.

Let us now describe the induction functor $u_{n*}$. Let $G\in\widehat{\cat{SM}}_{[n]}$
(recall that $\widehat{\cat{SM}}_{[n]}$ and $\widehat{\cat{SSM}_{n}}$
are equivalent; we shall now commit the crime of actually identifying
these two categories). Then\[
u_{n*}(G)(M\map N)=\Hom_{\widehat{\cat{SM}}_{[n]}}((\Pi T)^{n}(M\map N),G).\]
In particular, if $G$ is representable as a presheaf on $\cat{SM}$,
i.e.~if $G\in\cat{SM}_{[n]}$, then\[
u_{n*}(G)(M\map N)=\Hom_{\cat{SM}_{[n]}}((\Pi T)^{n}(M\map N),G),\]
or equivalently,\[
u_{n*}(G)(M\map N)=\Hom(C^{\infty}(G),\Omega_{[n]}(M\map N))\]
($\HOM(\rodd{n},\rodd{n})$-equivariant morphisms of algebras). In
particular for $n=1$,\[
u_{1*}(G)(M\map N)=\Hom(C^{\infty}(G),\Omega(M\map N))\]
(morphisms of differential graded algebras).

Presheaves of the form $u_{n*}(G)$, where $G\in\cat{SM}_{[n]}$,
will be called \emph{$n$-representable}. The category of $n$-representable
presheaves on $\cat{SSM}$ it thus equivalent to $\cat{SM}_{[n]}$. 

When we restrict a presheaf $F$ from $\cat{SSM}$ to $\cat{SSM}_{n}$
and induce it back to $\cat{SSM}$, we denote the result $\mathrm{app}_{n}F$;
in other words, $\mathrm{app}_{n}=u_{n*}\circ u_{n}^{*}$. We have
\[
(\mathrm{app}_{n}F)(M\map N)=\Hom_{\widehat{\cat{SM}}}((\Pi T)^{n}(M\map N),n\mbox{-jet of }F),\]
and if the $n$-jet of $F$ is representable,\[
(\mathrm{app}_{n}F)(M\map N)=\Hom(\Omega_{[n]}(F),\Omega_{[n]}(M\map N)).\]

If $m\geq n$ then \[
\mathrm{app}_{n}\circ\mathrm{app}_{m}\simeq\mathrm{app}_{m}\circ\mathrm{app}_{n}\simeq\mathrm{app}_{n}\]
(the isomorphisms are natural). The reason for this is Lemma \ref{lem:trivequiv}:
if we define auxiliary categories $\cat{SSM}_{\leq n}$ with object
$\rodd{k}\times N\map N$, $k\leq n$, then the pair $\cat{SSM}_{\leq n}\subset\cat{SSM}_{n}$
satisfies the assumptions of this lemma. We can thus substitute $\cat{SSM}_{n}$
with $\cat{SSM}_{\leq n}$ everywhere in this section. The categories
$\cat{SSM}_{\leq n}$ form an increasing chain, so our claim follows.

The morphism $\mathrm{id}_{\cat{SSM}}\map\mathrm{app}_{n}$ (coming
from the fact that $\mathrm{app}_{n}$ is the composition of an adjoint
pair) gives us a morphism\[
F\map\mathrm{app}_{n}F\]
for every $F\in\widehat{\cat{SSM}}$. Since the categories $\cat{SSM}_{\leq n}$
form an increasing chain, we get a chain\[
\mathrm{app}_{0}F\leftarrow\mathrm{app}_{1}F\leftarrow\mathrm{app}_{2}F\leftarrow\dots\]
which together with the morphisms $F\map\mathrm{app}_{n}F$ forms
a commutative diagram.

\section{Why simplicial manifolds?}

We could produce a simpler version of approximations etc.~(as described
above) by using finite sets in place of $\rodd{n}$. We would take
the full subcategory of $\cat{SSM}$ with objects of the form $S\times N\map N$,
where $S$ is a finite set and the arrow is the canonical projection.
Presheaves on this subcategory are functors\[
\cat{FS}^{o}\map\widehat{\cat{SM}},\]
where $\cat{FS}$ is the category of finite sets.

Simplicial (super)manifolds are functors\[
\Delta^{o}\map\cat{SM},\]
where $\Delta$ is the category of finite ordered sets. It would seem
more natural to use functors $\cat{FS}^{o}\map\cat{SM}$ for the purposes
of this paper. Indeed, the only reason why we used simplicial manifolds
is that they are better known.

\section{Examples \label{sec:Examples}}

\begin{example}
As the most trivial example, let us take the presheaf represented
by an object $(M\map N)\in\cat{SSM}$. This presheaf is $n$-representable
for any $n\geq1$, the corresponding object in $\cat{SM}_{[n]}$ is
$(\Pi T)^{n}(M\map N)$, hence $\Omega_{[n]}(F)=\Omega_{[n]}(M\map N)$.
The presheaf is $0$-representable only if $N$ is a point (which
is, however, the most interesting case).
\end{example}
{}

\begin{example}
\label{exa:group}This is the example of Section \ref{sec:Introductory-example}.
Let $G$ be a Lie group (or supergroup). Let \[
F(M\map N)=\{ G\mbox{-descent data on }M\map N\}.\]
In other words this is $F_{K_{\bullet}}$ for $K_{\bullet}$ the nerve
of $G$.

The 1-jet of $F$ is represented by $\mathfrak{g}[1]=\HOM(\rodd{1},G)/G$.
Differential forms on $F$, i.e.~the differential graded algebra
of functions on $\mathfrak{g}[1]$, is the Chevalley-Eilenberg complex
of $\mathfrak{g}$. We have $\mathrm{app}_{1}(F)(M\map N)=$ leafwise
flat $\mathfrak{g}$-connections on $M$. All higher approximations
of $F$ are isomorphic to $\mathrm{app}_{1}(F)$.
\end{example}
{}

\begin{example}
As a generalization of the previous example, let $G$ be a Lie groupoid
(or supergroupoid). We define \[
F(M\map N)=\{\textrm{Lie groupoid morphisms }(M\times_{N}M\rightrightarrows M)\map G\},\]
i.e.~$F$ is $F_{K_{\bullet}}$ for $K_{\bullet}$ the nerve of $G$.

The 1-jet of $F$ is represented by $A[1]$, where $A$ is the Lie
algebroid of $G$. The algebra of differential forms on $F$ is thus
$\Gamma(\bigwedge A)$, the Chevalley-Eilenberg complex of $A$. The
first approximation of $F$ is given by \[
\mathrm{app}_{1}(F)(M\map N)=\{\textrm{Lie algebroid morphisms }T(M\map N)\map A\}.\]
Higher approximations of $F$ are again isomorphic to $\mathrm{app}_{1}(F)$.
\end{example}
{}

\begin{example}
\label{exa:K2}Let us now consider $F=F_{K_{\bullet}}$ for $K_{\bullet}=K(U(1),2)$.
In other words, $F(M\map N)$ is the set of all maps \[
h:M\times_{N}M\times_{N}M\map U(1)\]
 such that\[
h(x,x,y)=h(x,y,y)=1\]
and\[
h(x,y,z)h(x,z,w)=h(x,y,w)h(y,z,w)\]
for any $x,y,z,w\in M$ lying over the same point of $N$. (In other
words, $h$ is descent data for $U(1)$-gerbes.)

The 1-jet of this $F$ is represented by $\mathbb{R}[2]$ (i.e.~the
corresponding $L_{\infty}$-algebra is $\mathbb{R}[1]$) and $\Omega(F)=\mathbb{R}[t]$
where $dt=0$ and $\deg t=2$. The first approximation of $F$ is
\[
\mathrm{app}_{1}(F)(M\map N)=\mathcal{Z}^{2}(M\map N)\]
(where $\mathcal{Z}^{2}$ denotes closed 2-forms). The morphism $F\map\mathrm{app}_{1}(F)$,
i.e.~the computation of a closed fibrewise 2-form $\omega$ out of
$h$, is given by \[
\omega(u,v)(x)=\mathcal{L}_{u}^{y}\mathcal{L}_{v}^{z}u(x,y,z)\mbox{ at }x=y=z.\]
Higher approximations of $F$ are not isomorphic to $\mathrm{app}_{1}(F)$
and they don't stabilize; they compute higher jets of $h$'s.

Similar example can be obtained from $K(U(1),n)$, where one gets
closed $n$-forms in place of closed 2-forms.
\end{example}
{}

\begin{example}
This example is a common generalization of Examples \ref{exa:group}
and \ref{exa:K2}. Recall that \emph{a Lie crossed module} is a pair
of Lie groups $G$ and $H$, with a morphism $m:H\map G$ and with
a morphism $\mu:G\map\mathit{Aut}(H)$ compatible with conjugations
on $G$ and $H$. Out of any Lie crossed module we can make a simplicial
manifold (its nerve) and the corresponding $F_{K_{\bullet}}$ looks
as follows: Elements of $F_{K_{\bullet}}(M\map N)$ are pairs of maps\[
g:M\times_{N}M\map G,\quad h:M\times_{N}M\times_{N}M\map H\]
such that\[
g(x,x)=1\]
\[
h(x,x,y)=h(x,y,y)=1\]
\[
g(x,y)\ g(y,z)=m(h(x,y,z))\ g(x,z)\]
\[
h(x,y,z)h(x,z,w)=h(x,y,w)\ \mu(g(x,y))(h(y,z,w)).\]
The 1-jet of $F_{K_{\bullet}}$ is represented by the differential
graded manifold $\mathfrak{f}[1]$, where $\mathfrak{f}$ is the following
DGLA: \[
\mathfrak{f}_{-1}=\mathfrak{h},\ \mathfrak{f}_{0}=\mathfrak{g},\ \mathfrak{f}_{i}=0\mbox{ otherwise,}\]
 the differential $\mathfrak{f}_{-1}\map\mathfrak{f}_{0}$ is $\mathfrak{h}\map\mathfrak{g}$
coming from $m$, the bracket on $\mathfrak{f}_{0}$ is that on $\mathfrak{g}$
and the bracket between $\mathfrak{f}_{0}$ and $\mathfrak{f}_{-1}$
is given by $\mu$. The functor $\mathrm{app}_{1}F$ is thus given
by\[
\mathrm{app}_{1}F(M\map N)=\{\mbox{flat fibrewise }\mathfrak{f}\mbox{-connections on }M\}.\]

\end{example}
{}

\begin{example}
Let $G$ be a Lie group acting by automorphisms on a commutative Lie
group $H$ and let $\phi:G^{n}\map H$ be a smooth group $n$-cocycle
for some $n\geq2$. The action $\mu:G\map\mathit{Aut}(H)$ and the
cocycle $\phi$ give us a simplicial manifold $K_{\bullet}$ which
is a fibre bundle over the nerve of $G$, with typical fibre $K(H,n-1)$.
The 1-jet of $F_{K_{\bullet}}$ is represented by the differential
graded manifold $\mathfrak{f}[1]$, where $\mathfrak{f}$ is the following
$L_{\infty}$-algebra:\[
\mathfrak{f}_{-n+2}=\mathfrak{h},\ \mathfrak{f}_{0}=\mathfrak{g},\ \mathfrak{f}_{i}=0\mbox{ otherwise,}\]
differential is zero, the bracket in $\mathfrak{f}_{0}$ is that of
$\mathfrak{g}$, the bracket between $\mathfrak{f}_{0}$ and $\mathfrak{f}_{-n+2}$
comes from $\mu$, and the only non-zero multiple bracket is\[
{\textstyle \bigwedge}^{n-1}\mathfrak{g}\map\mathfrak{h},\]
which is the image of $\phi$ under the Van Est map. The functor $\mathrm{app}_{1}F$
is thus given by flat fibrewise $\mathfrak{f}$-connections.
\end{example}
{}

Let us now pass to examples that are not given by a simplicial manifold.

\begin{example}
This is an extremely simple example, but it might be enlightening.
Let $F$ be given by \[
F(M\map N)=\mathcal{Z}^{k}(M\map N).\]
Then the 1-jet of $F$ is represented by the (1-dimensional) vector
superspace of closed $k$-forms on $\rodd{1}$, i.e.~by $\mathbb{R}[k]$.
This $F$ is $1$-representable.
\end{example}
{}

\begin{example}
This is an example of a presheaf on which a group acts. Let $G$ be
a Lie group, $\mathfrak{g}$ its Lie algebra, and let $F$ be given
by\[
F(M\map N)=\{\textrm{fibrewise }\mathfrak{g}\textrm{-connections on }M\}.\]
 On this $F$ acts the group (represented by) $G$, namely on $F(M\map N)$
acts the group $\Hom(M,G)$ by gauge transformations. The presheaf
$F$ is $1$-representable, the algebra of differential forms on $F$
is the Weil algebra $W(\mathfrak{g})$. Moreover, the action of $\HOM(\rodd{1},G)$
on differential forms gives the standard $G$-differential structure
on $W(\mathfrak{g})$.

If now $K$ is a manifold (or supermanifold) with $G$-action, the
algebra of differential forms on $F\times K$ is $W(\mathfrak{g})\otimes\Omega(K)$,
and the algebra of $\HOM(\rodd{1},G)$-invariant differential forms
is the basic subcomplex (whose cohomology is the equivariant cohomology
of $K$ if $G$ is compact).
\end{example}
{}

\begin{example}
Let us fix a supermanifold $Y$ and let $F$ be given by $F(M\map N)=\Hom(M\times_{N}M,Y)$.
Then the $n$-jet of $F$ is represented by $\HOM(\rodd{n}\times\rodd{n},Y)$,
and \[
\mathrm{app}_{n}(F)(M\map N)=\Gamma(j^{n}(M\map N,Y))\]
where $j^{n}(X,Y)\map X$ is the bundle of $n$-jets of maps $X\map Y$.
If for example $Y=\mathbb{R}$, the differential forms on $F$ form
the free graded commutative algebra generated by $x,\xi,\tau,t$ with
$\deg x=0$, $\deg\xi=\deg\tau=1$, $\deg t=2$, with differential
given by $dx=\xi$, $d\tau=t$.
\end{example}
{}

\begin{example}
Given $m\in\mathbb{N}$ we shall consider the presheaf given by

\[
F(M\map N)=\Gamma((T^{*})^{\otimes m}(M\map N))\]
and similar presheaves given by natural subbundles of $(T^{*})^{\otimes m}$
(it is convenient to decompose $(T^{*})^{\otimes m}$ with respect
to the action of the symmetric group $S_{m}$).

Let $\lambda$ be a Young diagram with $m$ squares, e.g.$$\includegraphics{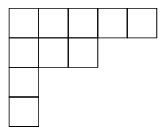}$$
for $m=10$, and let $W_{\left(\lambda\right)}$ be the corresponding
irreducible representation of $S_{m}$. If $V$ is a vector space,
let \[
V_{\lambda}=\Hom_{S_{m}}(W_{\left(\lambda\right)},V^{\otimes m}),\]
where $S_{m}$ acts on $V^{\otimes m}$ by permutations of factors
(recall that $V_{\lambda}$ is an irreducible representation of $GL(V)$).
This way we get a functor from the category of vector spaces to itself,
given by $V\mapsto V_{\lambda}$, and thus a presheaf \[
F_{\lambda}(M\map N)=\Gamma(T_{\lambda}^{*}(M\map N))\]

For any $\lambda$ and any $n$ the $n$-jet of $F_{\lambda}$ is
represented by the vector superspace $\Gamma(T_{\lambda}^{*}\rodd{n})$.
Let $c$ be the number of columns of $\lambda$. The presheaf $F_{\lambda}$
is $n$-representable for $n>c$, while $\mathrm{app}_{n}(F_{\lambda})(M\map N)=\{0\}$
for $n<c$. For $n=c$ we have \[
\mathrm{app}_{c}(F_{\lambda})(M\map N)=\Gamma(\widetilde{T_{\lambda}^{*}}(M\map N))\]
where $\widetilde{T_{\lambda}^{*}}$ is certain natural bundle containing
$T_{\lambda}^{*}$.

If $c=1$ (the case of differential forms) then $\widetilde{T_{\lambda}^{*}}=T_{\lambda}^{*}$;
the same is true when $\lambda$ has only one row (the case of symmetric
tensors). When $c=2$ we have the following result: $\widetilde{T_{\lambda}^{*}}$
has a natural increasing filtration\[
T_{\lambda}^{*}=A_{1}\subset A_{2}\subset\cdots\subset\widetilde{T_{\lambda}^{*}},\]
such that each $A_{i+1}/A_{i}$ is isomorphic to some $T_{\mu}^{*}$.
The rule for obtaining these $\mu$'s from $\lambda$ should be clear
from this picture:

$$\includegraphics{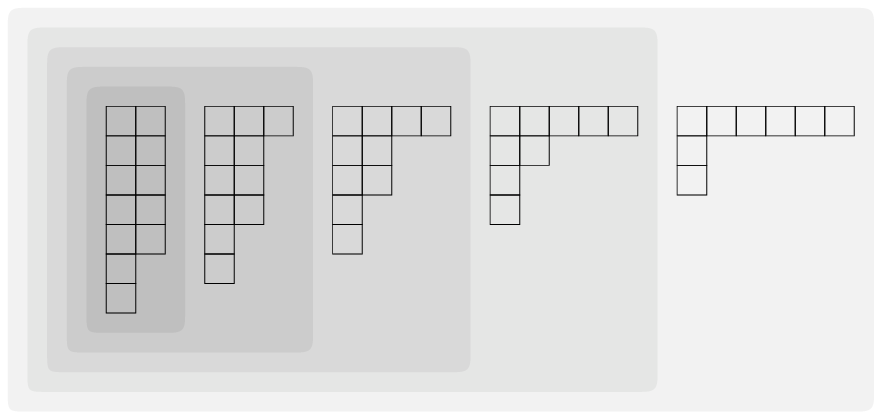}$$In other words, we keep removing squares
from the two column and adding them to the first row, till the second
column contains only one square.

This result enables us to decompose $\Omega_{[2]}(M)$ to indecomposable
representation of the super-semigroup $\HOM(\rodd{2},\rodd{2})$.
The irreducible (left) representations of $\HOM(\rodd{2},\rodd{2})$
are duals of $\Gamma(T_{\lambda}^{*}\rodd{2})$ for 2-column $\lambda$'s
(these are called generic irreducibles), and duals of the  spaces
of closed $k$-forms on $\rodd{2}$, $k\in\mathbb{N}$ (non-generic
irreducibles). The representation theory of $\HOM(\rodd{2},\rodd{2})$
is quite simple \cite{lei}: generic irreducibles never appear in
the composition series of reducible indecomposable representations,
and the only reducible indecomposable cyclic representations are duals
to the spaces of differential $k$-forms on $\rodd{2}$, $k\in\mathbb{N}$.
It gives the following decomposition of polynomial functions on $\HOM(\rodd{2},\rodd{2})$:\[
\Omega_{[2]}(\rodd{2})\cong\left(\bigoplus_{\lambda}\Gamma(T_{\lambda}^{*}\rodd{2})^{*}\otimes\Gamma(T_{\lambda}^{*}\rodd{2})\right)\oplus\frac{\bigoplus_{k}\Omega^{k}(\rodd{2})^{*}\otimes\Omega^{k}(\rodd{2})}{d(\bigoplus_{k}\Omega^{k}(\rodd{2})^{*}\otimes\Omega^{k-1}(\rodd{2}))}\]
(where the sum is over all 2-column $\lambda$'s) and thus for any
supermanifold $X$\[
\Omega_{[2]}(X)\cong\left(\bigoplus_{\lambda}\Gamma(T_{\lambda}^{*}\rodd{2})^{*}\otimes\Gamma(\widetilde{T_{\lambda}^{*}}X)\right)\oplus\frac{\bigoplus_{k}\Omega^{k}(\rodd{2})^{*}\otimes\Omega^{k}(X)}{d(\bigoplus_{k}\Omega^{k}(\rodd{2})^{*}\otimes\Omega^{k-1}(X))}.\]

Unfortunately, the representation theory of $\HOM(\rodd{n},\rodd{n})$
is wild for $n\geq3$ \cite{sh}, and we do not know how to decompose
the space of polynomial functions on this semigroup. We also do not
know the composition series of $\widetilde{T_{\lambda}^{*}}$ for
general $\lambda$ with more than 2 columns. Just as an example, here
is the composition series for some random (and quite simple) $\lambda$
with 3 columns:

$$\includegraphics{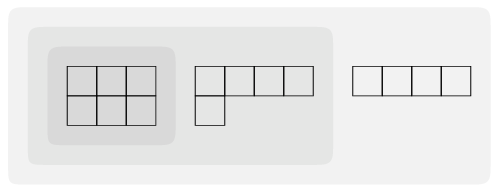}$$
\end{example}
\appendix

\section{\label{sec:app}Proof of Proposition \ref{pro:Kan}}

We use the definition of Kan simplicial manifolds from \cite{hen}.
Let $\Delta$ denote the category of finite ordered sets; a simplicial
set is a presheaf on $\Delta$. For any $n\in\mathbb{N}$ let $[n]$
denote the ordered set $0<1<\dots<n$, and also (as is usual) the
presheaf on $\Delta$ it represents. We use the standard notation
$F_{n}:=F([n])$ for any $F\in\widehat{\Delta}$. 

For any $k$, $0\leq k\leq n$, let $[n,k]\in\widehat{\Delta}$ (the
$k$-th horn of the simplex $[n]$) be the subpresheaf of $[n]$ defined
by\[
[n,k]_{m}=\{ f\in[n]_{m}=\Hom([m],[n]);\{0,\dots,k-1,k+1,\dots n\}\not\subset\mathrm{Im}(f)\}\]
(geometrically, $[n,k]$ is the union of the faces of the $n$-dimensional
simplex that do not contain the $k$-th vertex).

The inclusion $[n,k]\map[n]$ gives us a map $X_{n}=\Hom([n],X_{\bullet})\map\Hom([n,k],X_{\bullet})$
for any simplicial set $X_{\bullet}$. $X_{\bullet}$ is a \emph{Kan
simplicial set} if the map is surjective for all $n$'s and $k's$.
$X_{\bullet}$ is \emph{$m$-truncated} if the map is a bijection
for all $n\geq m$. The set $\Hom([n,k],X_{\bullet})$ will be denoted
$X_{n,k}$; it can be seen as a subset of $(X_{n-1})^{n}$.

A simplicial (super)manifold $K_{\bullet}$, i.e.~a functor\[
\Delta^{o}\map\cat{SM},\]
is \emph{Kan} if the map $K_{n}\map K_{n,k}$ is a surjective submersion
for any $n,k$. This definition requires checking that $K_{n,k}=\Hom([n,k],K_{\bullet})$
is actually a (super)manifold; it is done recursively in $n$ (see
\cite{hen} for details).

Let us now pass to the proof of Prop.~\ref{pro:Kan}. Let us first
describe the simple basic idea on sets, after that we shall pass to
(super)manifolds. Let $S$ be a set with a chosen element $*\in S$,
and let $S_{\bullet}$ be the nerve of the pair groupoid $S\times S\rightrightarrows S$,
i.e.~$S_{n}=S^{n+1}$, or more precisely,\[
S_{n}=\{\mbox{maps }[n]=\{0,1,\dots,n\}\map S\}.\]
 Let $X_{\bullet}$ be an $m$-truncated Kan simplicial set. Let us
inductively construct all morphisms $S_{\bullet}\map X_{\bullet}$
(basically by filling horns in all possible way). For any $k\in\mathbb{N}$
let us define a simplicial set $S_{\bullet}^{(k)}\subset S_{\bullet}$
by \[
S_{n}^{(k)}=\{\mbox{maps }f:[n]\map S\mbox{ s.t. }|\mathrm{Im}(f)|\leq k+1\mbox{ and }f(0)=*\mbox{ if }|\mathrm{Im}(f)|=k+1\};\]
they form a chain\[
S_{\bullet}^{(0)}\subset S_{\bullet}^{(1)}\subset S_{\bullet}^{(2)}\subset\dots\subset S_{\bullet}.\]
Let us denote \[
G^{(k)}:=\Hom(S_{\bullet}^{(k)},X_{\bullet}).\]
We shall notice that $G^{(k+1)}\map G^{(k)}$ is a surjection for
each $k$ and an bijection for $k\geq m$ (recall that $X$ is supposed
to be $m$-truncated) and that \[
\Hom(S_{\bullet},X_{\bullet})\map G^{(k)}\]
is also a bijection for $k\geq m$. Moreover, we shall describe how
to construct $G^{(k+1)}$ out of $G^{(k)}$.

First of all notice that $S_{\bullet}^{(k)}$ has non-degenerate simplices
only in dimensions $\leq k$; therefore any morphism $S_{\bullet}^{(k)}\map X_{\bullet}$
is uniquely determined by its dimension-$k$ part $\{*\}\times S^{k}=S_{k}^{(k)}\map X_{k}$.
The set $G^{(k)}$ is thus a subset of $\Hom(S^{k},X_{k})$. $\Hom(S_{\bullet},X_{\bullet})$
will therefore be constructed as a subset of $\Hom(S^{m},X_{m})$.

Obviously, $G^{(0)}\simeq X_{0}$. Suppose we already know $G^{(k)}$
and let us describe $G^{(k+1)}$. For any $g\in G^{(k)}$, $g:S^{k}\map X_{k}$,
we shall describe all $\tilde{g}$'s in $G^{(k+1)}$ lying over $g$.
For any $(s_{1},\dots,s_{k+1})\in S^{k+1}$ the $k+1$ elements $g(s_{1},\dots,s_{i-1},s_{i+1},\dots,s_{k+1})$
($1\leq i\leq k+1$) of $X_{k}$ form an element of $X_{k+1,0}$.
If $s_{k+1}\neq s_{k}$, let us choose an arbitrary element of $X_{k+1}$
lying over this element of $X_{k+1,0}$. If $s_{k+1}=s_{k}$, let
us take the image of $g(s_{1},\dots,s_{k})\in X_{k}$ under the $k+1$-th
degeneracy map $X_{k}\map X_{k+1}$. This way we choose an element
of $X_{k+1}$ for each $n+1$-tuple $(s_{1},\dots,s_{n+1})$, and
declare it to be $\tilde{g}(s_{1},\dots,s_{k+1})$. For $k\geq m$
the map $\tilde{g}$ is determined uniquely by $g$, since $X_{\bullet}$
is $m$-truncated.

Let us now pass from sets to supermanifolds. We shall need a little
lemma:

\begin{lem}
\label{lem:aux}Let $A=\rodd{k}$, $B=\rodd{\ell}$, and let $P$,
Q be any supermanifolds. Let $B\map A$ be an embedding and $P\map Q$
a surjective submersion, and let\[
\begin{array}{ccc}
\HOM(A,P) & \map & \HOM(A,Q)\\
\downarrow &  & \downarrow\\
\HOM(B,P) & \map & \HOM(B,Q)\end{array}\]
be induced by these two maps. Then the map \[
\HOM(A,P)\map\HOM(B,P)\times_{\HOM(B,Q)}\HOM(A,Q)\]
is a surjective submersion.
\end{lem}
The proof of the lemma is straightforward: if $P$, $Q$ are vector
superspaces and the map $P\map Q$ a linear map then one easily checks
that the result is a surjective linear map. The general (non-linear)
case then follows by the definition of a submersion.

Let us finally construct the supermanifold representing the $n$-jet
of $F_{K_{\bullet}}$, as a subsupermanifold of $\HOM((\rodd{n})^{m},K_{m})$.
We shall construct a chain of supermanifolds\[
H^{(k)}\subset\HOM((\rodd{n})^{k},K_{k})\]
analogous to $G^{(k)}$'s above, together with surjective submersions
\[
\dots\map H^{(2)}\map H^{(1)}\map H^{(0)},\]
where $H^{(k+1)}\map H^{(k)}$ is an isomorphism if $k\geq m$. $H^{(m)}$
will be the supermanifold we are looking for. 

The construction goes exactly as in the case of $G^{(k)}$'s; the
only difference is that we have to apply Lemma \ref{lem:aux}, namely
to\[
\begin{array}{ccc}
\HOM((\rodd{n})^{k+1},K_{k+1}) & \map & \HOM((\rodd{n})^{k+1},K_{k+1,0})\\
\downarrow &  & \downarrow\\
\HOM((\rodd{n})^{k},K_{k+1}) & \map & \HOM((\rodd{n})^{k},K_{k+1,0})\end{array},\]
where the canonical map $K_{k+1}\map K_{k+1,0}$ is a surjective submersion
by assuption, and the embedding $(\rodd{n})^{k}\map(\rodd{n})^{k+1}$
is defined to be the identity on the first $k-1$ $\rodd{n}$'s, and
to be the diagonal embedding of the last $\rodd{n}$ in $(\rodd{n})^{k}$
into the product of the two last $\rodd{n}$'s in $(\rodd{n})^{k+1}$.

$H^{(0)}$ is $K_{0}$. When we already have $H^{(k)}$, we have an
embedding \[
H^{(k)}\map\HOM((\rodd{n})^{k+1},K_{k+1,0})\]
and another embedding (given by the $k+1$-th degeneration map $K_{k}\map K_{k+1}$)
\[
H^{(k)}\map\HOM((\rodd{n})^{k},K_{k+1}).\]
It is easy to see that these two combine to an embedding\[
H^{(k)}\map\HOM((\rodd{n})^{k+1},K_{k+1,0})\times_{\HOM((\rodd{n})^{k},K_{k+1,0})}\HOM((\rodd{n})^{k},K_{k+1}).\]
Finally, $H^{(k+1)}$ is the preimage of this embedded $H^{(k)}$
under the submersion\[
\HOM((\rodd{n})^{k+1},K_{k+1})\map\HOM((\rodd{n})^{k+1},K_{k+1,0})\times_{\HOM((\rodd{n})^{k},K_{k+1,0})}\HOM((\rodd{n})^{k},K_{k+1}).\]

\end{document}